\documentclass[reqno, 12pt]{amsart}
\usepackage{amsmath}

\long\def\symbolfootnote[#1]#2{\begingroup%
\def\thefootnote{\fnsymbol{footnote}}\footnote[#1]{#2}\endgroup}

\newcommand\e{\epsilon}
\newcommand\R{{\mathbb R}}

\newtheorem{teo}{Theorem}[section]

\newtheorem{remark}[teo]{Remark}

\newtheorem{theorem}[teo]{Theorem}
\newtheorem{lemma}[teo]{Lemma}

\catcode`\@=11

\@addtoreset{equation}{section}
\@addtoreset{teo}{section}
\usepackage{amsbsy}
\def\RR{{\mathbb R}}

\newlength{\defbaselineskip}
\newcommand{\setlinespacing}[1]%
           {\setlength{\baselineskip}{#1 \defbaselineskip}}

\begin{document}
\title{Solutions of an elliptic system with a nearly critical exponent}
\thanks{This research was supported by
        FONDECYT 3040059}

\author{I.A. GUERRA}
\address{ Departamento de Matematica y C. C.,
Universidad  de Santiago,
 Ca\-si\-lla 307, Co\-rreo 2, Santiago, Chile}
 \email{\tt iguerra@usach.cl}

\subjclass{Primary,-----, Secondary,-----} \keywords{Semilinear elliptic system, asymptotic behaviour, critical exponent.  }


\begin{abstract}
Consider the problem
\begin{eqnarray*}
-\Delta u_\e &=& v_\e^p\,\,\quad  v_\e>0\quad \mbox{in}\quad \Omega,\\
-\Delta v_\e &=& u_\e^{q_\e}\quad  u_\e>0\quad \mbox{in}\quad \Omega, \\
u_\e&=&v_\e\:\:=\:\:0 \quad \mbox{on}\quad \partial \Omega,
\end{eqnarray*}
where $\Omega$ is a bounded convex domain in $\R^N,$ $N>2,$ with
smooth boundary $\partial \Omega.$ Here $p,q_\e>0,$ and
\begin{equation*}
\epsilon:=\frac{N}{p+1}+\frac{N}{q_\e+1}-(N-2).
\end{equation*}
This problem has positive solutions for $\e>0$ (with $pq_\e>1$)
and no non-trivial solution for $\e\leq 0.$ We study the
asymptotic behaviour of \emph{least energy} solutions as $\e\to
0^+.$ These solutions are shown to blow-up at exactly one point,
and the location of this point is characterized. In addition, the
shape and exact rates for blowing up are given.
\smallskip

\noindent {\sc R\'esum\'e.} Consid\'er\'e le probl\`{e}me
\begin{eqnarray*}
-\Delta u_\e &=& v_\e^p\,\,\quad  v_\e>0\quad \mbox{en}\quad \Omega,\\
-\Delta v_\e &=& u_\e^{q_\e}\quad  u_\e>0\quad \mbox{en}\quad \Omega, \\
u_\e&=&v_\e\:\:=\:\:0 \quad \mbox{sur}\quad \partial \Omega,
\end{eqnarray*}
o\`{u} $\Omega$ est un domaine convexe et born\'e de $\R^N,$
$N>2,$ avec la fronti\`{e}re r\'eguli\`{e}re $\partial \Omega.$
Ici $p,q_\e>0,$ et
\begin{equation*}
\epsilon:=\frac{N}{p+1}+\frac{N}{q_\e+1}-(N-2).
\end{equation*}
Ce probl\`{e}me a les solutions positives pour $\e>0$ (avec
$pq_\e>1$) et non pas de solution non-trivial pour $\e\leq 0.$
Nous \'etudions le comportement asymptotique de solutions
d'\emph{\'energie minimale}  quand $\e\to 0^+.$ Ces solutions
explosent en un seul point, et la localisation de ce point est
characteris\'e. De plus, la forme et le rythme d'explosion sont
donn\'es.

\end{abstract}

\maketitle

\section{Introduction}
We consider the elliptic system
\begin{eqnarray}\label{eq1}
-\Delta u_\e &=& v_\e^p\,\,\quad  v_\e>0\quad \mbox{in}\quad \Omega,\\
-\Delta v_\e &=& u_\e^{q_\e}\quad  u_\e>0\quad \mbox{in}\quad \Omega, \\
u_\e&=&v_\e\:\:=\:\:0 \quad \mbox{on}\quad \partial
\Omega,\label{eq2}
\end{eqnarray}
where $\Omega$ is a bounded convex domain in $\R^N,$ $N>2,$ with
smooth boundary $\partial \Omega.$ Here $p,q_\e>0,$ and
\begin{equation}\label{pareps}
\epsilon:=\frac{N}{p+1}+\frac{N}{q_\e+1}-(N-2).
\end{equation}
When $\epsilon\leq 0,$ there is no solution for
\eqref{eq1}-\eqref{eq2}, see \cite{M} and \cite{V}. On the other
hand when $\e>0,$ we can prove existence of  solutions  obtained
by the variational method. In fact, for $\e>0,$ the embedding
$W^{2,\frac{p+1}{p}}(\Omega)\hookrightarrow L^{q_\e+1}(\Omega)$ is
compact for any $q_\e+1>(p+1)/p,$ that is $pq_\e>1.$ Using this,
it is not difficult to show that there exists a function $\bar
u_\e$ positive solution of the variational problem
\begin{eqnarray}\label{Se}
S_\e(\Omega)=\inf
\left\{
\|\Delta u\|_{L^{\frac{p+1}{p}}(\Omega)} \:\mid\:u\in
W^{2,\frac{p+1}{p}}(\Omega),\quad
\|u\|_{L^{q_\e+1}(\Omega)}=1\right\},
\end{eqnarray}
see for example \cite{W}. This solution satisfies $-\Delta \bar
u_\e= \bar v_\e^p,\:\:-\Delta \bar v_\e=S_\e(\Omega) \bar
u_\e^{q_\e},$ in $\Omega$ and $\bar u_\e=\bar v_\e=0$  on
$\partial \Omega.$ After a suitable multiples of $\bar u_\e$ and
$\bar v_\e$, we obtain $u_\e$ and $v_\e$ solving
\eqref{eq1}-\eqref{eq2}.
We call $(u_\e,v_\e)$  the {\it least energy solution} to
\eqref{eq1}-\eqref{eq2}.
For others existence results, we
refer to \cite{CFM}, \cite{FF}, \cite{FR}, \cite{HvV93}, and
\cite{PvV}.

Note that by setting $v_\e=(-\Delta u_\e)^{1/p},$ we can write the
system \eqref{eq1}-\eqref{eq2} only in terms of $u_\e,$ that is
\begin{eqnarray}\label{eq3}
-\Delta (-\Delta u_\e)^{1/p} &=& u_\e^{q_\e}\quad  u_\e>0\quad \mbox{in}\quad \Omega \\
u_\e&=&\Delta u_\e=0 \quad \mbox{on}\quad \partial
\Omega.\label{eq4}
\end{eqnarray}

Concerning the least energy solutions, in \cite{W} it was proved
that $S_\e(\Omega)\to S$ as $\e\downarrow 0,$ where $S$ is
independent of $\Omega$ and moreover is the best Sobolev constant
for the inequality
\begin{equation}\label{ines}\|u\|_{L^{q+1}(\R^N)}\leq S^{-\frac p{p+1}}\|\Delta
u\|_{L^{\frac{p+1}{p}}(\R^N)}\end{equation} with $p,q,N$
satisfying
\begin{equation}\label{sobhyperbola}
\frac{N}{p+1}+\frac{N}{q+1}-(N-2)=0.
\end{equation}
This shows that the sequence $\{u_\e\}_{\e>0}$ of
least energy solutions of \eqref{eq3}-\eqref{eq4} satisfy  
\begin{eqnarray}\label{limmin}
S_\e(\Omega)=\frac{\int_\Omega |\Delta
u_\e|^{\frac{p+1}{p}}\,dx}{\|u_\e\|_{L^{q_\e+1}(\Omega)}^{\frac{p+1}{p}}}=S+o(1)\quad
\mbox{as}\quad \e\to 0.
\end{eqnarray}
Relation \eqref{sobhyperbola} defines a curve in $\R^2_+,$ for the
variables $p$ and $q.$ This curve is the so-called {\em Sobolev
Critical Hyperbola.}  By symmetry, we assume without restriction
that
\begin{eqnarray}\label{pinf}
2/(N-2)<p\leq p^*:=(N+2)/(N-2).
\end{eqnarray}
For each fixed value of $p$, the strict inequality gives a lower
bound for the dimension, i.e. $N>\max\{2,2(p+1)/p\}.$


In this article, we shall study in detail the asymptotic behaviour
of the variational solution $u_\e,$ of \eqref{eq3}--\eqref{eq4} as
$\e\downarrow 0$, that is, as $q_\e$ approaches from
\underline{below} to $q$ given by the Sobolev Critical Hyperbola
\eqref{sobhyperbola}.

The asymptotic behaviour of the equation \eqref{eq3}-\eqref{eq4}
as $\e\downarrow 0$ has already been studied for $p=p^*$ and
$p=1.$ Next we recall some of these results to introduce ours.

The case $p=p^*$ is equivalent to consider the single equation
$$
-\Delta u_\e= u_\e^{p^*-\e}\quad \mbox{in}\quad \Omega,\qquad
\mbox{and}\quad u_\e=0\quad \mbox{on}\quad \partial\Omega.
$$
This problem was studied in \cite{AP,FW,H,R}. There, exact rates
of blow-up were given and the location of blow-up points were
characterized. One key ingredient was the Pohozaev identity and
the observation that the solution $u_\e,$ scaled in the form
$\|u_\e\|_{L^{\infty}(\Omega)}^{-1}u_\e$ converges to $U$ solution
of
\begin{eqnarray}\label{eq1Csola}
-\Delta U&=& U^{p^*},\quad U(y)>0\quad\mbox{for}\quad y\in \RR^N\\
U(0)&=&1,\quad
 U\to 0,\quad\mbox{as}\quad |y|\to \infty,
 \label{eq2Csola}
\end{eqnarray}
which is unique, explicit, and radially symmetric. For the
location of blow-up and the shape of the solution away of the
singularity, it was proved that a scaled $u_\e,$ given by
$\|u_\e\|_{L^{\infty}(\Omega)}u_\e,$ converges to the Green's
function $G,$ solution of $-\Delta G(x,\cdot) = \delta_{x}$ in
$\Omega,$ $G(x,\cdot)=0$ on $\partial \Omega.$ The location of
blowing-up points are the critical points of $\phi(x):=g(x,x)$ (in
fact their minima, see  \cite{FW}), where $g(x,y)$ is the regular
part of $G(x,y)$, i.e
$$
g(x,y)=G(x,y)-\frac 1{(N-2)\sigma_N|x-y|^{N-2}}.
$$
In \cite{CG}, a similar result was proven in the case $p=1,$
$(N>4),$ where the problem is reduced to study
\eqref{eq1Csola}--\eqref{eq2Csola} with the operator $\Delta^2$
instead of $-\Delta$. Both cases give the blow-up rate
$$
\epsilon \|u_\epsilon\|_{L^\infty(\Omega)}^2\to C\quad
\mbox{as}\quad\epsilon\to 0.
$$
for some explicit $C:=C(p,N,\Omega)>0.$ We can ask ourselves if
this behaviour is universal. We will see later that this is only a
coincide.

Mimicking the above argument, we will study the asymptotic
behaviour of the solution $u_\e$ of \eqref{eq3}--\eqref{eq4} as
$\e\downarrow 0.$ We shall show that
$\|u_\e\|_{L^{\infty}(\Omega)}^{-1}u_\e$ converges, as
$\e\downarrow 0,$ to the solution $U$ of the problem
\begin{eqnarray}\label{eq1C}
-\Delta U&=& V^p,\quad V(y)>0\quad\mbox{for}\quad y\in \RR^N\\
-\Delta V&=& U^{q},\quad U(y)>0\quad\mbox{for}\quad y\in \RR^N \\
U(0)&=&1,\quad
 U\to 0,\quad V\to 0\quad\mbox{as}\quad |y|\to \infty.
 \label{eq2C}
\end{eqnarray}
In \cite{CLO}, it was proved that $U$ and $V$ are radially
symmetric, if $p\geq 1$ and $U\in L^{q+1}(\R^N)$ and $V\in
L^{p+1}(\R^N)$. This is the case when considering least energy
solutions, see details in section \ref{preliminaries}. Thus
$U(r):=U(y)$ and $V(r):=V(y)$ with $r=|y|,$ moreover $U$ and $V$
are unique, and decreasing in $r,$ see \cite{HvV,W}.  There exist
no explicit form of $(U,V)$ for all $p\geq 1,$ however to carry
out the analysis it is sufficient to know the asymptotic behaviour
of $(U,V)$ as $r\to\infty,$ which was studied in \cite{HvV}. They
found
\begin{eqnarray}
&&\quad\lim\limits_{r\to \infty}r^{N-2}V(r)=a
\quad\mbox{and}\quad
\label{hulshof2} \begin{cases}
 \lim\limits_{r\to \infty}r^{N-2}U(r)=b\quad \mbox{if}\quad p>\frac{N}{N-2} \cr
 \lim\limits_{r\to \infty}\frac{r^{N-2}}{\log r}U(r)=b\quad \mbox{if}\quad p=\frac{N}{N-2} \cr
 \lim\limits_{r\to \infty}r^{p(N-2)-2}U(r)=b\quad \mbox{if}\quad \frac
 2{N-2}<p<\frac{N}{N-2}.
\end{cases}
\end{eqnarray}

The aim of this paper is to show the following results.

\begin{theorem}\label{tmain1}
Let $u_\e$ be a least energy solution  of \eqref{eq3}--\eqref{eq4}
and $p\geq 1.$ Then
%

a) there exists $x_0\in \Omega$ such that,
after passing to a subsequence, we have
$$
{\rm i)}\quad u_\e\to 0 \in C^1(\Omega\setminus\{x_0\}),
\qquad {\rm ii)}\quad v_\e=|\Delta u_\e|^{\frac 1{p}}\to 0 \in
C^1(\Omega\setminus\{x_0\})
$$
as $\e\to 0$ and
$$
{\rm iii)}\quad |\Delta u_\e|^{\frac{p+1}{p}}\to \|V\|_{L^{p+1}(\RR^N)}^{p+1}\delta_{x_0}\quad \mbox{as}\quad \e\to 0
$$
in the sense of distributions.

b) $x_0$ is a critical point of  \begin{eqnarray}
\phi(x)&:=&g(x,x)\quad\mbox{if}\quad p\in
[N/(N-2),(N+2)/(N-2))\quad \mbox{and} \\\tilde \phi(x)&:=&\tilde
g(x,x)\quad\mbox{if}\quad p\in (2/(N-2),N/(N-2))\end{eqnarray} for
$x\in \Omega.$
The function $\tilde g(x,y)$
is defined for $p\in (2/(N-2),N/(N-2))$ by
$$
\tilde g(x,y)= \tilde G(x,y)-\frac 1{(p(N-2)-2)(N-p(N-2))(N-2)^p\sigma_N^p|x-y|^{p(N-2)-2}}
$$
where 
$-\Delta \tilde G(x,\cdot) = G^p(x,\cdot)$  in $\Omega,$ $\tilde
G(x,\cdot)=0$ on $\partial \Omega.$
\end{theorem}
We observe that regularity of  $\tilde \phi$ is needed to compute
its critical points in b).
We show next that $\tilde \phi$ is regular.
By definition of $\tilde G$, we have
\begin{eqnarray}\label{limtilde}
\lim\limits_ {y\to x}|x-y|^{(p-1)(N-2)}\Delta\tilde
g(x,y)=-\frac{pg(x,x)}{((N-2)\sigma_N)^{p-1}}
\end{eqnarray}
for $x\in \Omega.$ Thus $-\Delta\tilde g(x,\cdot)\in
L^{q}(\Omega)$ for any $q\in (N/2, N/(p(N-2)-N+2)).$ This implies,
by regularity, that $\tilde g(x,\cdot)\in L^\infty(\Omega)$ and
therefore $\tilde \phi(x)=\tilde g(x,x),$ $x\in \Omega$ is
bounded. In addition, we define
\begin{eqnarray}\label{ghat}
\hat g(x,y)= \tilde
g(x,y)+\frac{pg(x,x)|x-y|^{N-p(N-2)}}{(N-p(N-2))(2N-p(N-2)-2)((N-2)\sigma_N)^{p-1}}
\end{eqnarray}
and we have for any $x\in \Omega$ that
\begin{eqnarray}\label{limtilde2}
\lim\limits_ {y\to x}|x-y|^{(p-2)(N-2)}\Delta\hat
g(x,y)=-\frac{p(p-1)g(x,x)}{((N-2)\sigma_N)^{p-2}}.
\end{eqnarray}
Thus $\hat g(x,y)$ is regular in $y$ for $x$ fixed. Since
$N>p(N-2),$ we take first $y=x$ in \eqref{ghat} and then the
gradient and we find $\nabla_x\tilde g(x,x)=\nabla_x\hat g(x,x).$
Hence $\tilde \phi(x)$ is regular.

To state the next theorems we denote
$$
\alpha=\frac{N}{q+1}\quad \mbox{and}\quad \beta=\frac{N}{p+1},
$$
so the critical hyperbola \eqref{sobhyperbola} takes the form
$\alpha+\beta=N-2.$

\begin{theorem}\label{tmain2} Let the assumptions of Theorem \ref{tmain1} be satisfied.
Then
$$
\begin{cases}
\quad \lim\limits_{\epsilon\to 0^+}\epsilon
\|u_\epsilon\|_{L^\infty(\Omega)}^\frac{(N-2)}{\alpha}=
S^{\frac{1-pq}{p(q+1)}}
\|U\|_{L^q(\RR^N)}^{q}\|V\|_{L^p(\RR^N)}^{p}|\phi(x_0)|\quad
\mbox{if}\quad p>\frac{N}{N-2} \cr \quad \lim\limits_{\epsilon\to
0^+}\epsilon
\frac{\|u_\epsilon\|_{L^\infty(\Omega)}^\frac{(N-2)}{\alpha}}{\log(\|u_\epsilon\|_{L^\infty(\Omega)})}=\frac{1}{\alpha}a^{\frac{N}{N-2}}
S^{\frac{1-pq}{p(q+1)}}\|U\|_{L^q(\RR^N)}^{q}|\phi(x_0)| \quad
\mbox{if}\quad p=\frac{N}{N-2} \cr \quad \lim\limits_{\epsilon\to
0^+}\epsilon
\|u_\epsilon\|_{L^\infty(\Omega)}^\frac{p(N-2)-2}{\alpha}=
S^{\frac{1-pq}{p(q+1)}}\|U\|_{L^q(\RR^N)}^{q(p+1)}|\tilde
\phi(x_0)|\quad \mbox{if}\quad p<\frac{N}{N-2}.
\end{cases}
$$
\end{theorem}

In particular taking  $p=p^*,$ and using (\ref{sobhyperbola}) we
find that $q=p^*.$ We recover the results in \cite{H,R}, that is
\begin{equation}\label{universal}
 \epsilon \|u_\epsilon\|_{L^\infty(\Omega)}^2\to C\quad
\mbox{as}\quad\epsilon\to 0,
\end{equation} for some explicitly
given $C>0$. See also \cite{AP} for the case $\Omega=B_R.$

When $N>4,$ we can take $p=1,$ and use (\ref{pinf}) to find that
$q=(N+4)/(N-4).$ Here we recover the result in \cite{BE,CG}, where
they prove that \eqref{universal} holds for some $C>0.$

\begin{theorem}\label{tmain3}
Let the assumptions of Theorem \ref{tmain1} be satisfied. Then
\begin{eqnarray}
&&\quad \lim\limits_{\e\to 0}\|u_\epsilon\|_{L^\infty(\Omega)}v_\epsilon(x)=\|U\|_{L^q(\RR^N)}^qG(x,x_0), \quad\mbox{and}\\
&&\begin{cases} \quad \lim\limits_{\e\to 0}\|u_\epsilon\|_{L^\infty(\Omega)}^{\frac
\beta{\alpha}}u_\epsilon(x)=\|V\|_{L^p(\RR^N)}^pG(x,x_0)\quad \mbox{if}\quad p>\frac{N}{N-2} \cr \quad  \lim\limits_{\e\to
0}\frac{\|u_\epsilon\|_{L^\infty(\Omega)}^{\frac \beta{\alpha}}}{\log \|u_\epsilon\|_{L^\infty(\Omega)}}u_\epsilon(x)=
\frac 1{\alpha}a^{\frac N{N-2}} G(x,x_0) \quad \mbox{if}\quad
p=\frac{N}{N-2} \cr \quad \lim\limits_{\e\to
0}\|u_\epsilon\|_{L^\infty(\Omega)}^{\frac
1{\alpha}(\beta+p(N-2)-N)}
u_\epsilon(x)=\|U\|_{L^q(\RR^N)}^{pq}\tilde G(x,x_0)\quad
\mbox{if}\quad p<\frac{N}{N-2}
\end{cases}
\label{c1a2}
\end{eqnarray} where all the convergences in $C^{1,\alpha}(w)$ with $w$ any neighborhood of $\partial \Omega$ not containing $x_0.$
\end{theorem}

\begin{remark}
For $p<\frac{N}{N-2},$ the convergence in \eqref{c1a2} can be
improved to $C^{3,\alpha}(\omega).$ See the proof of the theorem.
\end{remark}

\begin{remark}
By \eqref{convcontinuos}, we find that $\lim\limits_{\e\to
0}\|v_\e\|_{L^\infty(\Omega)}=V(0)\lim\limits_{\e\to 0}
\|u_\e\|_{{L^\infty(\Omega)}}^{\frac \beta{\alpha}}.$ So, in
addition, when $p=1$ we have that
$$
\epsilon \|v_\epsilon\|_{L^\infty(\Omega)}^{2(N-4)/N}\to
CV(0)^{2(N-4)/N}\quad \mbox{as}\quad\epsilon\to 0.
$$
\end{remark}

We can extend these results to the problem
\begin{eqnarray}\label{eq5}
-\Delta (-\Delta u_\e)^{1/p} &=& u_\e^{q}+\e u_\e\quad  u_\e>0\quad \mbox{in}\quad \Omega \\
u_\e&=&\Delta u_\e=0 \quad \mbox{on}\quad \partial \Omega\label{eq6}
\end{eqnarray}
with $\epsilon\to 0.$ The existence of positive solutions for this
problem can be found in \cite{HvV93} and \cite{PvV} in the case of
a ball. See \cite{Ham} for related results for $p=1.$

\begin{theorem}
Let the assumptions of Theorem \ref{tmain1} be satisfied.
Then
\begin{align}
&\lim\limits_{\epsilon\to 0^+}\epsilon
\|u_\epsilon\|_{L^\infty(\Omega)}^{2-\frac{2}{\alpha}}=
\|U\|^{-2}_{L^2(\RR^N)}
\|U\|_{L^q(\RR^N)}^{q}\|V\|_{L^p(\RR^N)}^{p}|\phi(x_0)|\quad
\mbox{if}\quad p>\frac{N}{N-2}\nonumber\\ &\qquad\qquad
\qquad\qquad \qquad\qquad \qquad\qquad \mbox{and}\quad \alpha >
1,\:N>4\\ &\lim\limits_{\epsilon\to 0^+}\epsilon \frac
{\|u_\epsilon\|_{L^\infty(\Omega)}^{2-\frac{2}{\alpha}}}{\log(\|u_\epsilon\|_{L^\infty(\Omega)})}
=\frac{1}{\alpha}a^{\frac{N}{N-2}}\|U\|^{-2}_{L^2(\RR^N)}
\|U\|_{L^q(\RR^N)}^{q}|\phi(x_0)| \quad \mbox{if}\quad
p=\frac{N}{N-2} \nonumber\\ &\qquad\qquad \qquad\qquad
\label{align2} \qquad\qquad \qquad\qquad \mbox{and}\quad
2-\frac{2}{\alpha}=3-(\frac{N}{N-2})^2>0,\\
& \lim\limits_{\epsilon\to 0^+}\epsilon
\|u_\epsilon\|_{L^\infty(\Omega)}^{2-\frac{2+N-p(N-2)}{\alpha} }=
\|U\|^{-2}_{L^2(\RR^N)}\|U\|_{L^q(\RR^N)}^{q(p+1)}|\tilde
\phi(x_0)|\quad \mbox{if}\quad
\frac{N+4}{2(N-2)}<p<\frac{N}{N-2}\quad \nonumber\\ &\qquad\qquad
\qquad\qquad \qquad\qquad \qquad\qquad \mbox{and}\quad \alpha>
\frac{2+N-p(N-2)}{2} \\
&\lim\limits_{\epsilon\to 0^+}\epsilon
\log(\|u_\epsilon\|_{L^\infty(\Omega)})=
\frac{\|U\|_{L^q(\RR^N)}^{q}\|V\|_{L^p(\RR^N)}^{p}}{b^2}|\phi(x_0)|\quad
\mbox{if}\quad N=4,\:p=q=3,
\end{align}
\begin{eqnarray}
\lim\limits_{\epsilon\to 0^+}\epsilon \log
(\|u_\epsilon\|_{L^\infty(\Omega)})\|u_\epsilon\|_{L^\infty(\Omega)}^{\frac{3-q}{2}}&=&
\frac{\|U\|_{L^q(\RR^N)}^{q(p+1)}}{b^2}|\tilde \phi(x_0)|\quad
\mbox{if}\quad p=\frac{N+4}{2(N-2)}<\frac{N}{N-2},\nonumber \\
\label{align5} \qquad\qquad &&\qquad\qquad \qquad\qquad
\mbox{and}\quad q\leq 3
\end{eqnarray}
\end{theorem}
Note that $N>4$ (integer) is equivalent to $3-(N/(N-2))^2>0$ and
also to $(N+4)/(2(N-2))<N/(N-2).$ This implies that \eqref{align2}
holds for $p=N/(N-2)$ and $N>4,$ and  \eqref{align5} holds for
$p=\frac{N+4}{2(N-2)},$ $q\leq 3$ and $N>4.$


For example, $p=1$ gives $q+1=2N/(N-4)$ and provided that $N>8,$
we get
$$
\lim\limits_{\epsilon\to 0^+}\epsilon
\|u_\epsilon\|_{L^\infty(\Omega)}^{\frac{2(N-8)}{N-4}}= C_1|\tilde
\phi(x_0)|.
$$
For $N=8$ and  $p=1,$ we have
$$ \lim\limits_{\epsilon\to
0^+}\epsilon \log (\|u_\epsilon\|_{L^\infty(\Omega)})= C_1|\tilde
\phi(x_0)|.
$$

\section{Preliminaries}\label{preliminaries}
Before proving the main theorem, we need some properties of
$u_\e.$ Using that $u_\e$ is a minimizing sequence, we have
$$
\int\limits_\Omega(\Delta u_\e)^{\frac{p+1}{p}}\,dx=
\int\limits_\Omega v_\e\Delta u_\e\,dx=\int\limits_\Omega
u_\e\Delta v_\e\,dx=\int\limits_\Omega u_\e^{q_\e+1}\,dx.
$$
Then $
[S+o(1)]\|u_\e\|^{\frac{p+1}{p}}_{L^{q_\e+1}(\Omega)}=\|u_\e\|_{L^{q_\e+1}(\Omega)}^{q_\e+1}
$ implies
\begin{eqnarray}\label{convqe}
\lim\limits_{\e\to 0}\int\limits_\Omega
u_\e^{q_\e+1}\,dx=S^{\frac{pq-1}{p(q+1)}}.
\end{eqnarray}

\begin{lemma}
The minimizing sequence  $u_\e$ of \eqref{limmin} is such that
$$
\|u_\e\|_{L^{\infty}(\Omega)}\to\infty
$$
moreover $\|(-\Delta u_\e)^{1/p}\|_{L^{\infty}(\Omega)}=\|v_\e\|_{L^{\infty}(\Omega)}\to\infty\quad \mbox{as}\quad \e\to 0.
$
\end{lemma}
\begin{proof}
 If $\|u_\e\|_{L^{\infty}(\Omega)}\to\infty$ then by regularity, we find $\|v_\e\|_{L^{\infty}(\Omega)}\to\infty$,
 see \cite[Theorem 3.7]{GT}.
Now, assume that $ \|u_\e\|_{L^{\infty}(\Omega)}\leq M$ and
$\|v_\e\|_{L^{\infty}(\Omega)}\leq M,$ by elliptic regularity, we
have that
$$
  \|v_\e\|_{C^{2+\alpha}(\bar \Omega)}\leq M
  \quad
  \mbox{and}\quad
  \|u_\e\|_{C^{2+\alpha}(\bar \Omega)}\leq M
$$
with $\alpha\in(0,1)$ and some constant $M.$ This implies that
there exists $u^*,v^*\in  C^2(\bar \Omega),$ such that
$$
u_\e\to u^*\quad\mbox{in}\quad C^2(\bar \Omega) ,\quad v_\e\to
v^*\quad\mbox{in}\quad C^2(\bar \Omega)\quad\mbox{as}\quad \e\to
0.
$$
Hence $u^*$ satifies
$$
0\neq \int\limits_\Omega(\Delta
u^*)^{\frac{p+1}{p}}\,dx=S\left[\int\limits_\Omega
(u^*)^{q+1}\,dx\right]^{\frac{(p+1)}{p(q+1)}}
$$
which contradicts that $S$ can be achieved by a
 minimizer in a bounded domain, see
\cite{W}. In other words there exists no non trivial solution for
\begin{eqnarray}\label{eq1crit}
-\Delta u^*&=& (v^*)^p,\quad  v>0\quad \mbox{in}\quad \Omega\\
-\Delta v^* &=& (u^*)^{q},\quad  u>0\quad \mbox{in}\quad \Omega \\
 u^*&=&v^*=0 \quad \mbox{on}\quad \partial \Omega \label{eq2crit}
\end{eqnarray}
in a convex bounded domain, with $p,q$ satisfying (\ref{sobhyperbola}), see \cite{M},\cite{V}.
\end{proof}

For any $\e>0,$ let $(u_\e,v_\e)$ be a solution of (\ref{eq1}--\ref{eq2}).
By the Pohozaev inequality, see \cite{M} or \cite{V}, we have for any $\tilde \alpha,\tilde \beta\in\RR$ that
\begin{eqnarray}
\left(\frac{N}{q_\e+1} -\tilde \alpha\right)\int\limits_\Omega u_\e^{q_\e+1}\,dx +\left(\frac N{p+1}-\tilde \beta\right)\int\limits_\Omega v_\e^{p+1}\,dx \\
+(N-2-\tilde \alpha-\tilde \beta)\int\limits _\Omega(\nabla u_\e,\nabla v_\e)\,dx =-\int\limits_{\partial \Omega}(\nabla
u_\e,n)(\nabla v_\e,x-y)\,ds.
\end{eqnarray}
We choose $\tilde \alpha+\tilde \beta=N-2$, $\tilde \alpha=\alpha$
and so $\tilde\beta=\beta.$ This implies that
\begin{eqnarray}\label{epoh}
\e\int\limits_\Omega u_\e^{q_\e+1}\,dx=-\int\limits_{\partial
\Omega}\frac{\partial  u_\e}{\partial n}\frac{\partial
v_\e}{\partial n}(n,x-y)\,ds.
\end{eqnarray}

Since $ u_\e$ becomes unbounded  as $\e\to 0$ we choose $\mu=\mu(\e)$ and  $x_\e\in\Omega$  such that $$ \mu^{\alpha_\e}
u_\epsilon(x_\e)=1
$$
where $\alpha_\e=N/(q_\e+1).$ Note that $\mu\to 0$ as $\e\to 0$.

First we claim that $x_\e$ stays away from the boundary. This is
consequence of moving plane method and interior estimates
\cite{DLN}, \cite{GNN}. Let $\phi_1$ the positive eigenvalue of
$(-\Delta,H_0^1(\Omega)),$ normalized to $\max\limits_{x\in
\Omega} \phi_1(x)=1.$ Since $p\geq 1,$ multiplying by $\phi_1$ we
obtain
\begin{eqnarray*}
\lambda_1\int\limits_\Omega u_\e\phi_1= \int\limits_\Omega v_\e^p\phi_1\geq 2 \lambda_1\int\limits_\Omega v_\e\phi_1-C\int\limits_\Omega \phi_1 \\
\lambda_1\int\limits_\Omega v_\e\phi_1= \int\limits_\Omega
u_\e^{q_\e}\phi_1\geq 2 \lambda_1\int\limits_\Omega
u_\e\phi_1-C\int\limits_\Omega \phi_1
\end{eqnarray*}
for some $C=C(p,q,\lambda_1)>0$. Hence $ \int\limits_\Omega
u_\e\phi_1\leq (C/\lambda_1)\int\limits_\Omega\phi_1 $ which
implies $\int_{\Omega'} u_\e\leq C(\Omega')$ with $\Omega'\subset
\Omega$ and $\int_{\Omega'} v_\e\leq C(\Omega').$ Using the moving
planes method \cite{GNN}, we find that there exist $t_0\alpha>0$
such that
$$
u_\e(x-t\nu)\quad\mbox{and}\quad  v_\e(x-t\nu)\quad\mbox{are
nondecreasing for}\quad  t\in [0,t_0],
$$
$\nu\in\RR^N$ with  $|\nu|=1,$ and $(\nu,n(x))\geq \alpha$ and $x\in \partial\Omega.$
Therefore we can find $\gamma,\delta$ such that for any $x\in \{z\in \bar \Omega\colon d(z,\partial\Omega)<\delta\,\}=\Omega_\delta$ there exists a measurable set $\Gamma_x$ with (i) $meas(\Gamma_x)\geq \gamma,$ (ii) $\Gamma_x\subset \Omega\setminus\bar \Omega_{\delta/2},$  and (iii) $u_\e(y)\geq u_\e(x)$ and $v_\e(y)\geq v_\e(x)$ for any $y\in \Gamma_x.$ Then for any $x\in \Omega_\delta,$ we have
\begin{eqnarray*}
u_\e(x)&\leq& \frac 1{meas(\Gamma_x)}\int\limits_{\Gamma_x}u_\e(y)dy\leq \frac 1{\gamma}\int\limits_{\Omega_\delta}u_\e\leq C(\Omega_\delta),\quad \mbox{and} \\
v_\e(x)&\leq& \frac 1{meas(\Gamma_x)}\int\limits_{\Gamma_x}v_\e(y)dy\leq \frac 1{\gamma}\int\limits_{\Omega_\delta}v_\e\leq
C(\Omega_\delta).
\end{eqnarray*}
Hence if $u_\e(x_\e)\to \infty,$ this implies that $x_\e$ will stay out of $\Omega_\delta$ a neighbordhood of the boundary. This proves the claim.


Let $x_\e\to x_0\in \Omega.$ We define a family of rescaled
functions
\begin{eqnarray}
u_{\e,\mu}(y)&=&\mu^{\alpha_\e} u_\epsilon(\mu^{1-\e/2} y+x_\e)\quad  \\
v_{\e,\mu}(y)&=&\mu^\beta  v_\epsilon(\mu^{1-\e/2}y+x_\e)
\end{eqnarray}
and find using the definitions of $\e,$ $\alpha_\e$ and $\beta,$ that
\begin{eqnarray}\label{emusystem1}
-\Delta u_{\e,\mu}&=&v_{\e,\mu}^p\mu^{\alpha_\e+2-\e-p\beta}= v_{\e,\mu}^p\quad\mbox{in}\quad \Omega_\e \\
-\Delta v_{\e,\mu}&=&u_{\e,\mu}^{q_\e}\mu^{\beta+2-\e-q_\e\alpha_\e}= u_{\e,\mu}^{q_\e}\quad\mbox{in}\quad \Omega_\e \\
u_{\e,\mu}&=& v_{\e,\mu}=0\quad \mbox{on}\quad \partial \Omega_\e. \label{emusystem2}
\end{eqnarray}
By equicontinuity and using Arzela-Ascoli, we have that
\begin{eqnarray}\label{convcontinuos}
u_{\e,\mu}\to U\quad \mbox{and}\quad v_{\epsilon,\mu}\to V\quad \mbox{as}\quad \epsilon\to 0.
\end{eqnarray}
in $C^{2}(K)$ for any $K$ compact in $\RR^N,$ where $(U,V)$
satisfies \eqref{eq1C}--\eqref{eq2C}.
 Now extending $u_{\e,\mu}$ and $v_{\e,\mu}$ by zero outside $\Omega_\e$ and using \eqref{convqe}, by  the argument
 in \cite{S} or \cite{W}, we have that $u_{\e,\mu}\to \bar U$ strongly (up to a subsequence) in
 $W^{2,\frac{p+1}{p}}(\R^N).$ In the limit $\bar U\in L^{q+1}(\R^N)$ and $\bar V:=(-\Delta \bar U)^{\frac 1{p}}\in L^{p+1}(\R^N),$
 and they satisfy \eqref{eq1C}--\eqref{eq2C}. Since $p\geq 1,$
the solution $(\bar U,\bar V)$ is unique and radially symmetric,
see \cite{CLO}. In addition the radial solutions are unique
\cite{HvV,W}, so $\bar U\equiv U$ and $\bar V\equiv V$,
consequently
\begin{eqnarray}\label{convueve}
\int\limits_{\RR^N}[u_{\e,\mu}-U]^{q+1}(y)\,dy \to 0
\quad
\int\limits_{\RR^N} [v_{\e,\mu}-V]^{p+1}(y)\,dy\to 0.
\end{eqnarray}

\begin{lemma}\label{deltalemma} There exists $\delta>0$ such that
 $$\delta\leq  \mu^\e\leq 1. $$
\end{lemma}
\begin{proof}
Since $\mu\to 0,$ we have $\mu^\e\leq 1.$ By \eqref{convueve}, we
get $ \int\limits_{B_1}u_{\e,\mu}^{q_\e+1}\,dx\geq M, $ but
\begin{eqnarray}
M\leq \int\limits_{B_1}u_{\e,\mu}^{q_\e+1}\,dx=\mu^{\e
N/2}\int\limits_{|y-x_\e|\leq \mu^{1-\e/2}}u_{\e}^{q_\e+1}(y)\,dy
\leq \mu^{\e N/2}\int\limits_{\Omega}u_{\e}^{q_\e+1}(y)\,dy
\end{eqnarray}
Using the convergence \eqref{convqe}, we obtain the result.
\end{proof}

\begin{lemma}\label{mainlemma}
There exists $K>0$ such that the solution $(u_{\e,\mu},v_{\e,\mu})$ satisfies
\begin{eqnarray}\label{boundgs}
u_{\e,\mu}(y)\leq KU(y)\quad v_{\e,\mu}(y)\leq KV(y) \quad \forall
y\in \RR^N.
\end{eqnarray}
\end{lemma}
We prove this lemma in section \ref{mainlemma}.


\begin{lemma} There exists a constant $C>0$ such that
\begin{eqnarray}\label{eboundmu}
\e\leq C\mu^{N-2}h(\mu) \quad \mbox{with}\quad h(\mu)=\begin{cases} 1 \quad \qquad\quad\mbox{for}\quad p>N/(N-2) \cr |\log(\mu)|
\quad\quad \mbox{for}\quad p=N/(N-2) \cr \mu^{(p(N-2)-N)} \quad \mbox{for}\quad p<N/(N-2).
\end{cases}
\end{eqnarray}
\end{lemma}

\begin{proof}
We will establish the following
$$
\int\limits_{\partial \Omega}\frac{\partial u_\e}{\partial n}\frac{\partial  v_\e}{\partial n}(n,x)\, dx\leq C \mu^{N-2}h(\mu)
$$
and from here the result follows applying \eqref{epoh}. We claim that
$$
\left|\frac{\partial u_\e}{\partial n}\right|\leq C\mu^{\alpha_\e}
\quad \left|\frac{\partial v_\e}{\partial n}\right|\leq
C\mu^{\beta }h(\mu)
$$
In the following $M$ is a positive constant that can vary from
line to line and we shall use systematically Lemma
\ref{deltalemma}.

For $p>N/(N-2),$ using that $-p\beta+N=\beta,$ we have
$$
\int\limits_\Omega v_\e^p(x)\,dx\leq M \mu^{-p\beta+N(1-\e/2)}\int\limits_{\R^N} V^p(y)\,dy\leq M \mu^{\beta }
$$
and by \eqref{boundgs} there exists $M>0$ such that
\begin{eqnarray}\label{cotave}
v_\e^p(x)\leq
M\frac{\mu^{\beta+p(N-2)-N-p(N-2)\e/2}}{|x-x_0|^{p(N-2)}}.
\end{eqnarray}
for $x\neq x_0.$ Using  that $\beta<\beta+p(N-2)-N,$ by Lemma \ref{Lpreg} we find  $|\frac{\partial v_\e}{\partial n}|\leq
C\mu^{\beta}$. For $u_\e,$ using that $-q_\e\alpha_\e+N=\alpha_\e,$
$$
\int\limits_\Omega u_\e^{q_\e}\,dx\leq M \mu^{-q_\e\alpha_\e+N(1-\e/2)}\int\limits _{\RR^N} U^q(y)\,dy\leq M \mu^{\alpha_\e }
$$
and  by \eqref{boundgs} there exist $M>0$ such
\begin{eqnarray}\label{uecota1}
u_\e^{q_e}(x)\leq M
\frac{\mu^{-q_\e\alpha_\e+q_\e(N-2)-q_\e(N-2)\e/2}}{|x-x_0|^{q_\e(N-2)}}
\end{eqnarray}
for $x\neq x_0.$ Using  that $\alpha_\e<\alpha_\e-N+q_\e(N-2),$ by Lemma \ref{Lpreg}, we obtain  $|\frac{\partial u_\e}{\partial
n}|\leq C\mu^{\alpha_\e}$.

For $p<N/(N-2),$ we have
\begin{eqnarray}
\int\limits_\Omega v_\e^p\,dx\leq M \mu^{-p\beta+p(N-2)(1-\e/2)}\lim\limits_{\mu\to 0} \frac
1{\mu^{(p(N-2)-N)(1-\e/2)}}\int\limits_{B_{\frac 1{\mu^{1-\e/2}}}(x_\e)}  V^p(y)\,dy
\\
\leq M \mu^{\beta +(p(N-2)-N) }
\end{eqnarray}
and for $x\neq x_0,$ we find \eqref{cotave} for $v_\e$
and for $u_\e$ we have
$$
\int\limits_\Omega u_\e^{q_\e}\leq M \mu^{-q_\e\alpha_\e+N(1-\e/2)}\int \limits_{\RR^N}  U^q(y)\,dy\leq M \mu^{\alpha_\e}
$$
and by \eqref{boundgs} there exist $M>0$ such that \begin{eqnarray}\label{uecota2}
u_\e^{q_e}(x)\leq M
\frac{\mu^{-q_\e\alpha_\e+q_\e(p(N-2)-2)-q_\e(p(N-2)-2)\e/2}}{|x-x_0|^{q_\e(p(N-2)-2)}}
\end{eqnarray}
for $x\neq x_0$. From these estimates we prove the claim applying Lemma \ref{Lpreg} and noting that
$\alpha_\e<\alpha_\e-N+q_\e(p(N-2)-2)+(p+1)\e/\alpha_\e.$. For the case $p=N/(N-2),$ we proceed as before noting that
$$
\int\limits_\Omega v_\e^p\,dx\leq M \mu^{-p\beta+N(1-\e/2)}|\log(\mu)|\lim\limits_{\mu\to 0} \frac
1{|\log(\mu)|}\int\limits_{B_{\frac 1{\mu^{1-\e/2}}}(x_\e)}  V^p(y)\,dy\leq M |\log(\mu)|\mu^{\beta}
$$
and for $x\neq x_0$ we have \eqref{cotave}. Similarly to
\eqref{uecota2}, we obtain that for $x\neq x_0$, there exist $M>0$
such that
\begin{eqnarray}\label{uecota3}
u_\e^{q_e}(x)\leq M
\frac{\mu^{-q_\e\alpha_\e+q_\e(N-2)-q_\e(N-2)\e/2}}{|x-x_0|^{q_\e(N-2)}}\log(|x-x_0|\mu^{-1+\e/2})^{q_\e}.
\end{eqnarray}
Using this and proceeding and before we prove the claim and  the lemma follows.
\end{proof}

\begin{lemma}\label{convmuto1}
$$|\mu^\e-1|=O(\mu^{N-2}h(\mu)\log\mu)$$
\end{lemma}
\begin{proof} By the theorem of the mean $|\mu^\e-1|=|\mu^{s\e}\e\log\mu|$ for some $s\in(0,1)$ and therefore \eqref{eboundmu} gives the result.
\end{proof}

\section{Proof of Theorem \ref{tmain2} and \ref{tmain3}}


\begin{proof}[Proof of Theorem \ref{tmain3}]  We start by proving the case $p>\frac{N}{N-2}.$  We have
\begin{eqnarray}
-\Delta (\|u_\e\|^{\frac{\beta}{\alpha}}_{L^{\infty}(\Omega)}u_\e)&=
& \|u_\e\|^{\frac{\beta}{\alpha}}_{L^{\infty}(\Omega)} v_\e^p \quad \mbox{in}\quad \Omega,\label{th15a}\\
-\Delta (\|u_\e\|_{L^{\infty}(\Omega)}v_\e)&=& \|u_\e\|_{L^{\infty}(\Omega)} u_\e^{q_\e} \quad \mbox{in}\quad \Omega,\\
u_\e&=&v_\e=0\quad \mbox{on}\quad \partial \Omega.
\end{eqnarray}
We integrate the right hand side of \eqref{th15a}
$$
\int\limits_{\Omega}\|u_\e\|^{\frac{\beta}{\alpha}}_{L^{\infty}(\Omega)}
v_\e^p\, dx=\mu^{-(p+1)\beta+N+N\e /2}\int\limits_{\Omega_\e}
v_{\e,\mu}^p(y)\,dy.
$$
But $ N-(p+1)\beta=0$, so using \eqref{boundgs} by dominated
convergence and Lemma \ref{convmuto1}, we get
$$
\lim\limits_{\e\to
0}\int\limits_{\Omega}\|u_\e\|^{\frac{\beta}{\alpha}}_{L^{\infty}(\Omega)}
v_\e^p\, dx= \int\limits_{\RR^N} V^p(y)
\,dy=\|V\|_{L^p(\R^N)}<\infty.
$$
Similarly, now using
\begin{eqnarray}\label{uelq}
\int\limits_{\Omega}\|u_\e\|_{L^{\infty}(\Omega)} u_\e^{q_\e}\,
dx= \mu^{-(q_\e+1)\alpha_\e+N+N\e /2} \int\limits_{\Omega_\e}
u_{\e,\mu}^{q_\e}\, dx \to \|U\|_{L^q(\R^N)}<\infty
\end{eqnarray}
as $\e\to 0.$ Also using the bound \eqref{boundgs}, we find
$$
\|u_\e\|^{\frac{\beta}{\alpha}}_{{L^\infty(\Omega)}} v_\e^p (x) \leq
\frac{M\mu^{-(p+1)\beta+p(N-2)-p(N-2)\e/2}}{|x-x_0|^{p(N-2)}}
$$
for $x\neq x_0$ and some $M>0.$  But $-(p+1)\beta+p(N-2)>0$ and Lemma \ref{deltalemma} then
$\|u_\e\|^{\frac{\beta}{\alpha}}_{L^\infty(\Omega)} v_\e^p(x)\to 0$ for $x\neq x_0.$  Also we have
$$
\|u_\e\|_{{L^\infty(\Omega)}} u_\e^{q_\e} (x) \leq \frac{M\mu^{-(q_\e+1)\alpha_\e+q_\e(N-2)-q_\e(N-2)\e/2}}{|x-x_0|^{q_\e(N-2)}}.
$$
for $x\neq x_0$ and some $M>0.$  But $-(q_\e+1)\alpha_\e+q_\e(N-2)>0$ and Lemma \ref{deltalemma} then
$\|u_\e\|_{L^\infty(\Omega)} u_\e^{q_\e}(x)\to 0$ for $x\neq x_0.$

From here we have
$$
-\Delta (\|u_\e\|^{\frac{\beta}{\alpha}}_{L^\infty(\Omega)} u_\e)
\to \|V\|_{L^p(\RR^N)}^{p}\delta_{x=x_0}\quad\mbox{and}\quad
-\Delta (\|u_\e\|_{L^\infty(\Omega)} v_\e) \to
\|U\|_{L^q(\RR^N)}^{q}\delta_{x=x_0}
$$
in the sense of distributions in $\Omega,$ as $\e\to 0.$ Let
$\omega$ be any neighborhood of $\partial \Omega$ not containing
$x_0.$ By regularity theory, see Lemma \ref{Lpreg}, we find
$$
\|\|u_\e\|^{\frac{\beta}{\alpha}}_{L^\infty(\Omega)} u_\e
\|_{C^{1,\alpha}(w)}\leq C\left
[\|\|u_\e\|^{\frac{\beta}{\alpha}}_{L^\infty(\Omega)} v_\e^p
\|_{L^1(\Omega)}+\|\|u_\e\|^{\frac{\beta}{\alpha}}_{L^\infty(\Omega)}
v_\e^p\|_{L^\infty(w)}\right]
$$
and a similar bound for $\|\|u_\e\|_{L^\infty(\Omega)} v_\e
\|_{C^{1,\alpha}(w)}.$  Consequently
\begin{eqnarray}\label{utog}
\|u_\e\|^{\frac{\beta}{\alpha}}_{L^\infty(\Omega)} u_\e\to
\|V\|_{L^p(\RR^N)}^{p}G\quad \mbox{in} \quad
C^{1,\alpha}(w)\quad\mbox{as}\quad \e\to 0.
\end{eqnarray}
and
\begin{eqnarray}\label{vtog}
\|u_\e\|_{L^\infty(\Omega)} v_\e\to \|U\|_{L^q(\RR^N)}^{q}G\quad
\mbox{in} \quad C^{1,\alpha}(w)\quad\mbox{as}\quad \e\to 0.
\end{eqnarray}
For the case $p<N/(N-2),$ we proceed as before and we have \eqref{uelq} and the bound
$$
\|u_\e\|_{{L^\infty(\Omega)}} u_\e^{q_\e} (x) \leq \frac{M\mu^{-(q_\e+1)\alpha_\e+q_\e(p(N-2)-2)-q_\e(p(N-2)-2)\e/2}}{|x-x_0|^{q_\e(p(N-2)-2)}}.
$$
for $x\neq x_0$ and some $M>0.$ Using  that $-(q_\e+1)\alpha_\e+q(p(N-2)-2)=2(p+1)>0$ and Lemma \ref{deltalemma},  we get
$\|u_\e\|^{\frac{\beta}{\alpha}}_{L^\infty(\Omega)} u_\e^{q_\e}(x)\to 0$ for $x\neq x_0$ and hence
\begin{eqnarray}\label{vtog2}
\|u_\e\|_{L^\infty(\Omega)} v_\e\to
\|U\|_{L^q(\RR^N)}^{q}G\quad\mbox{in} \quad C^{1,\alpha}(w)
\quad\mbox{as}\quad \e\to 0.
\end{eqnarray}
Now we claim that
\begin{eqnarray}\label{utog2}
\|u_\e\|^{\frac{1}{\alpha}(\beta+p(N-2)-N)}_{L^\infty(\Omega)}
u_\e\to \|U\|_{L^q(\RR^N)}^{pq}\tilde G\quad\mbox{in}\quad
C^{1,\alpha}(w)\quad\mbox{as}\quad \e\to 0.
\end{eqnarray}
We have
$$
-\Delta
(\|u_\e\|^{\frac{1}{\alpha}(\beta+p(N-2)-N)}_{L^\infty(\Omega)}
u_\e)=\|u_\e\|^{\frac{1}{\alpha}(\beta+p(N-2)-N)}_{L^\infty(\Omega)}v_\e^p=\|u_\e\|_{L^\infty(\Omega)}^pv_\e^p.
$$
Since the last term converges to $(\|U\|_{L^q(\RR^N)}^{q}G)^p$ in
$C^{1,\alpha}(\omega)$ as $\e\to 0$ and $p\geq 1,$ we have
$$
\|u_\e\|^{\frac{1}{\alpha}(\beta+p(N-2)-N)}_{L^\infty(\Omega)}
u_\e\to \|U\|_{L^q(\RR^N)}^{pq}\tilde G\quad\mbox{in}\quad
C^{3,\alpha}(w)\quad \mbox{as}\quad \e\to 0.$$

For the remaining case $p=N/(N-2),$ we have as $\e\to 0,$ the
convergence
$$
\int\limits_{\Omega}\frac{\|u_\e\|^{\frac{\beta}{\alpha}}_{L^{\infty}(\Omega)}}{|\log
(\|u_\e\|_{L^{\infty}(\Omega)})|}v_\e^p\, dx=
\frac{\mu^{-(p+1)\beta+N+N\e /2}}{\alpha_\e|\log
(\mu)|}\int\limits_{\Omega_\e} v_{\e,\mu}^{p}\, dy \to \frac
1{\alpha}\lim\limits_{r\to
\infty}V(r)^{\frac{N}{N-2}}r^N=\frac{a^{\frac{N}{N-2}}}{\alpha}.
$$
and the pointwise bound for $x\neq x_0$
$$
\frac{\|u_\e\|^{\frac{\beta}{\alpha}}_{{L^\infty(\Omega)}} }{|\log
(\|u_\e\|^{\frac{\beta}{\alpha}})|}v_\e^p (x) \leq
\frac{M\mu^{-p(N-2)\e/2}}{\log(\mu)|x-x_0|^{p(N-2)}}.
$$
By Lemma \ref{deltalemma},
$\frac{\|u_\e\|^{\frac{\beta}{\alpha}}_{{L^\infty(\Omega)}}
}{|\log (\|u_\e\|^{\frac{\beta}{\alpha}})|} v_\e^{p}(x)\to 0$ for
$x\neq x_0.$ Writing
$$
-\Delta \left(\frac{\|u_\e\|^{\frac{\beta}{\alpha}}_{L^{\infty}}}{|\log (\|u_\e\|^{\frac{\beta}{\alpha}})|}
u_\e\right)=\frac{\|u_\e\|^{\frac{\beta}{\alpha}}_{L^{\infty}}}{|\log (\|u_\e\|^{\frac{\beta}{\alpha}})|} v_\e^p,
$$
we observe that the last term converges to $\delta_{x=x_0}.$ By
Lemma \ref{Lpreg}, we have
 $$
 \frac{\|u_\e\|^{\frac{\beta}{\alpha}}_{L^{\infty}}}{|\log (\|u_\e\|^{\frac{\beta}{\alpha}})|} u_\e\to \frac{a^\frac N{N-2}}{\alpha}G\quad\mbox{in}\quad C^{1,\alpha}(w)\quad\mbox{as}\quad \e\to 0,
 $$
and clearly we have \eqref{vtog} using \eqref{uecota3}. This
completes the proof of the theorem.
\end{proof}

\begin{proof}[Proof of Theorem \ref{tmain2}]
For $p>N/(N-2)$ we have
\begin{eqnarray*}
\e\|u_\e\|^{\frac{N-2}{\alpha}}\int\limits_\Omega
u_\e^{q_\e+1}\,dx=\int\limits_{\partial
\Omega}(\|u_\e\|^{\frac{\beta}{\alpha}}_{{L^\infty(\Omega)}}\nabla
u_\e,n)(\|u_\e\|_{{L^\infty(\Omega)}}\nabla v_\e,n)(n,x-y)\,ds
\end{eqnarray*}
By \eqref{utog} and \eqref{vtog},
\begin{eqnarray*}
\lim\limits_{\e\to
0}\e\|u_\e\|^{\frac{N-2}{\alpha}}\int\limits_\Omega
u_\e^{q_\e+1}\,dx=\|V\|_{L^p(\RR^N)}^{p}\|U\|_{L^q(\RR^N)}^{q}\int\limits_{\partial
\Omega}\frac{\partial G(x,x_0)}{\partial  n}\frac{\partial
G(x,x_0)}{\partial n}(n,x-x_0)\,ds.
\end{eqnarray*}
Also for the case $p<N/(N-2),$ using
\begin{eqnarray*}
\e\|u_\e\|^{\frac{1}{\alpha}(p(N-2)-2)}_{L^\infty(\Omega)}\int\limits_\Omega
u_\e^{q_\e+1}\,dx
\qquad\qquad\qquad\qquad\qquad\qquad\qquad\qquad\\=\int\limits_{\partial
\Omega}(\|u_\e\|^{\frac{1}{\alpha}(\beta+p(N-2)-N)}_{{L^\infty(\Omega)}}\nabla
u_\e,n)(\|u_\e\|_{{L^\infty(\Omega)}}\nabla v_\e,n)(n,x-y)\,ds
\end{eqnarray*}
and \eqref{utog2} and \eqref{vtog2}, we get
\begin{eqnarray*}
\lim\limits_{\e\to
0}\e\|u_\e\|^{\frac{1}{\alpha}(p(N-2)-2)}_{L^\infty(\Omega)}\int\limits_\Omega
u_\e^{q_\e+1}\,dx
\qquad\qquad\qquad\qquad\qquad\qquad\\=\|U\|_{L^q(\RR^N)}^{q(p+1)}\int\limits_{\partial
\Omega}\frac{\partial \tilde G(x,x_0)}{\partial  n}\frac{\partial
G(x,x_0)}{\partial n}(n,x-x_0)\,ds.
\end{eqnarray*}
The case $p=N/(N-2)$ is analogous.

The proof of the theorems follows from the next lemma.
\end{proof}

\begin{lemma} We have the following identities
\begin{eqnarray*}
i) \int\limits_{\partial \Omega}\frac{\partial G(x,x_0)}{\partial n}\frac{\partial G(x,x_0)}{\partial n}(n,x-x_0)\,ds=-(N-2)g(x_0,x_0)
\end{eqnarray*}
and
\begin{eqnarray*}
ii) \int\limits_{\partial \Omega}\frac{\partial \tilde
G(x,x_0)}{\partial n}\frac{\partial G(x,x_0)}{\partial n}
(n,,x-x_0)\,ds=-\frac N{q+1}\tilde g(x_0,x_0)
\end{eqnarray*}
\end{lemma}
\begin{proof} i) was proven in \cite{BP}, see also \cite{H}.
To prove ii) we follow a similar procedure. From \cite{M,V}, for
any $y\in \R^N,$ we have the following identity
\begin{eqnarray*}
\int\limits_{\Omega'} \Delta u(x-y,\nabla v)+\Delta v(x-y,\nabla u) - (N-2)(\nabla u,\nabla v)dx=\\
\int_{\partial \Omega'} \frac{\partial u }{\partial n}(x-y,\nabla
v)+\frac{\partial v }{\partial n}(x-y,\nabla u)-(\nabla u,\nabla
v)(x-y,n)\, ds
\end{eqnarray*}
where $\Omega'=\Omega\setminus B_r$ with $r>0.$ For a system
$-\Delta v= 0$ and  $-\Delta u= v^p,$ in $\Omega',$ the identity
takes the form
\begin{eqnarray}\nonumber
\int\limits_{\Omega'} \frac{N}{p+1}v^{p+1}- \bar a v^{p+1} dx=\int\limits_{\partial \Omega'} \frac 1{p+1}v^{p+1}(x-y,n)\,ds\\
+\int\limits_{\partial \Omega'} \frac{\partial u }{\partial
n}\left[(x-y,\nabla v)+\bar a v \right]+\frac{\partial v
}{\partial n}\left[(x-y,\nabla u)+\bar b u\right]- (\nabla
u,\nabla v)(x-y,n)\,ds \label{Mident}
\end{eqnarray}
with $\bar a+\bar b=N-2.$ Let $y=0$, choose $\bar a=N/(p+1)$ and
take $v=G(x,0)$ and $u=\tilde G(x,0).$ Using that $u=v=0$ on
$\partial \Omega,$ and so $\nabla u= (\nabla u,n)n$ and $\nabla v=
(\nabla v,n)n$ on $\partial \Omega,$, we obtain
 \begin{eqnarray*}
 \int\limits_{\partial \Omega}\frac{\partial \tilde G }{\partial n}\frac{\partial \tilde G }{\partial n}(x,n)ds
&=&\int\limits_{\partial B_r} \frac 1{p+1}G^{p+1}(x,n)+\frac{\partial \tilde G }{\partial n}[(x,\nabla G)+\frac{N}{p+1} G]\, ds \\
&&+\int\limits_{\partial B_r} \frac{\partial G }{\partial
n}[(x,\nabla \tilde G)+\frac{N}{q+1}\tilde G]-(\nabla \tilde
G,\nabla G)(x,n)\, ds.
\end{eqnarray*}
Let $k=p(N-2)$ and $\Gamma=\sigma_N(N-2).$ For $|x|=r,$ we have
\begin{eqnarray*}
\nabla \tilde G =-\frac 1{\Gamma^p(N-k)}|x|^{-k}x+\nabla \tilde
g,\quad \nabla G =-\frac 1{\sigma_N}|x|^{-N}x+\nabla g,
\end{eqnarray*}
\begin{eqnarray*}
\frac{\partial \tilde G }{\partial n}=-\frac
1{\Gamma^p(N-k)}|x|^{1-k}+(\nabla \tilde g,n),\quad \frac{\partial
G }{\partial n}=-\frac 1{\sigma_N}|x|^{1-N}+(\nabla g,n)
\end{eqnarray*}
\begin{eqnarray*}
 (x,\nabla \tilde G)+\frac{N}{q+1}\tilde G&=&(\frac{N}{(q+1)(k-2)}-1)\frac 1{\Gamma^p(N-k)}|x|^{2-k}+(x,\nabla \tilde g)+\frac{N}{q+1}\tilde g \\
  (x,\nabla G)+\frac{N}{p+1}G&=&(\frac{N}{p+1}-(N-2))\frac 1{\Gamma}|x|^{2-N}+(x,\nabla g)+\frac{N}{p+1} g
\end{eqnarray*}
\begin{eqnarray*}
(\nabla \tilde G,\nabla G)=\frac{|x|^{-k-N+2}}{\sigma_N
\Gamma^p(N-k)}-\frac{(\nabla g,x)}{\Gamma^p(N-k)}|x|^{(2-N)p}
-\frac {(\nabla \tilde g,x)}{\sigma_N}|x|^{-N}+(\nabla \tilde g,\nabla g)
\end{eqnarray*}
and
$$
\frac 1{p+1}G^{p+1}=\frac 1{p+1}\left[\frac
{1}{\Gamma^p}|x|^{-k}-\Delta \tilde g\right]\left[\frac
1{\Gamma}|x|^{2-N}+g\right]
$$
From here, we check that terms with $|x|^{3-N-k}$ cancel out
others integral tends to 0 since the integrand are $o(|x|^{1-N})$
and only remain one term of order $|x|^{1-N},$ giving
 \begin{eqnarray*}
 \int\limits_{\partial \Omega}\frac{\partial \tilde G }{\partial n}\frac{\partial G }{\partial n}(x,n)\,ds =-\lim\limits_{r\to 0}\frac 1{\sigma_N r^{N-1}}\int\limits_{\partial B_r}\frac{N}{q+1}\tilde g ds= -\frac{N}{q+1}\tilde g(0,0).
\end{eqnarray*}
\end{proof}

\begin{proof}[Proof of Theorem \ref{tmain1}]  a) The part ii) follows from Theorem \ref{Lpreg},
$$
\||\Delta u_\e|^{\frac 1{p}}\|_{C^{1,\alpha}(\omega)}\leq \|u_\e^{q_\e}\|_{L^1(\Omega)}+\|u_\e^{q_\e}\|_{L^\infty(\omega)}.
$$
and estimates \eqref{uecota1}, \eqref{uecota2}, and
\eqref{uecota3}. Part i) follows from
$$
\|u_\e\|_{C^{1,\alpha}(\omega)}\leq \|v_\e^{p}\|_{L^1(\Omega)}+\|v_\e^{p}\|_{L^\infty(\omega)}.
$$
and estimate \eqref{cotave}. Finally iii) follows combining ii) with the convergence
$$
\int\limits_{\R^N}|\Delta u_\e|^{\frac{p+1}{p}}\,dx=\int\limits_{\R^N}v_\e^{p+1}\,dx\to \|V\|^{p+1}_{L^{p+1}(\RR^N)}.
$$
as $\e\to 0.$ This completes part a).

For part b), note that from  \eqref{epoh}, we have the vectorial equality $
\int\limits_{\partial \Omega}(\nabla u_\e,\nabla v_\e)n\,ds=0.
$
In the limit for $p\geq N/(N-2),$ we get
\begin{eqnarray}\label{vec1}
\int\limits_{\partial \Omega}(\nabla G(x,x_0),\nabla G(x,x_0))n\,ds=0
\end{eqnarray}
and similarly for $p< N/(N-2),$ we obtain
\begin{eqnarray}\label{vec2}
\int\limits_{\partial \Omega}(\nabla \tilde G(x,x_0),\nabla
G(x,x_0))n\,ds=0
\end{eqnarray}

But we have the following result.
\begin{lemma} For every $x_0\in \Omega$
\begin{eqnarray}\label{vec3}
\int\limits_{\partial \Omega}(\nabla G(x,x_0) ,n)(\nabla G(x,x_0),n) n\,ds=-\nabla \phi(x_0)
\end{eqnarray}
and
\begin{eqnarray}\label{vec4}
\int\limits_{\partial \Omega}(\nabla \tilde G(x,x_0) ,n)(\nabla
(\Delta \tilde G(x,x_0))^{1/p},n) n\,ds=-\nabla \tilde \phi(x_0).
\end{eqnarray}
\end{lemma}

Hence combining \eqref{vec1} with \eqref{vec3}, and \eqref{vec2}
with \eqref{vec4}, we complete the proof of part b) and the
theorem is proven.
\end{proof}

\begin{proof}[Proof of the Lemma.]
Equality \eqref{vec3} was proved in \cite{BP} and \cite{H}. To
prove  \eqref{vec4}, by \eqref{Mident} we have
\begin{eqnarray*}
 \int\limits_{\partial \Omega}\frac{\partial \tilde G }{\partial n}\frac{\partial G }{\partial
 n}n\,ds
=\int\limits_{\partial B_r}\left\{ \frac
1{p+1}G^{p+1}n+\frac{\partial \tilde G }{\partial n}\nabla G\, +
\frac{\partial G }{\partial n}\nabla \tilde G-(\nabla \tilde
G,\nabla G)n\,\right\} ds.
\end{eqnarray*}
Using $\int\limits_{\partial B_r} n=0,$ we get
\begin{eqnarray}\nonumber
 \int\limits_{\partial \Omega}\frac{\partial \tilde G }{\partial n}\frac{\partial G }{\partial
 n}n\,ds=\frac 1{(p+1)r^{N-1}}\int\limits_{\partial B_r}\left\{\frac 1{\Gamma^p}r^{N-k-1}g-\Delta \tilde g
\frac 1{\Gamma} r-\Delta \tilde g g r^{N-1}\right\}n\,ds \\
\nonumber +\frac 1{r^{N-1}}\int\limits_{\partial B_r}\{(\nabla
\tilde g,n)\nabla g+(\nabla g,n)\nabla \tilde g-(\nabla \tilde g,\nabla g)n\}r^{N-1}\,ds \\
-\frac 1{r^{N-1}}\int\limits_{\partial
B_r}\left\{\frac{1}{\sigma_N}\nabla \tilde
g+\frac{r^{N-k}}{\Gamma^p(N-k)}\nabla g\right\}\,ds. \label{Grto0}
\end{eqnarray}
We use the regular $\hat g(x,0)$ instead of $\tilde g(x,0).$ Thus

\begin{eqnarray}\label{gradghat}
\nabla \hat g(x,0) &=& \nabla \tilde g(x,0)
+\frac{pg(0,0)}{\Gamma^{p-1}(2N-k-2)}|x|^{N-k-2}x, \\ \Delta \hat
g(x,0) &=& \Delta \tilde g(x,0)
+\frac{pg(0,0)}{\Gamma^{p-1}}|x|^{N-k-2}.\label{lapghat}
\end{eqnarray}
But $ g(x,0)=g(0,0)+(\nabla g(0,0),x)+o(|x|^2) $ and
$$
\int\limits_{\partial
B_r}r^{-k}g(x,0)n\,ds=\int\limits_{\partial
B_1}r^{N-k-1}g(0,0)n\,ds+\int\limits_{\partial B_1}r^{N-k}(\nabla
g(0,0),y) n\,ds+ o(r^{N-k+1})
$$
where $y=x/r.$ Clearly the first integral in the r.h.s is zero and
the other terms tends to zero as $r\to 0.$ Hence
\begin{eqnarray}\label{limghat}
\lim\limits_{r\to 0}\frac 1{r^{N-1}} \int\limits_{\partial
B_r}r^{N-k-1}g(x,0)n\,ds=0.
\end{eqnarray}
We replace \eqref{gradghat} and \eqref{lapghat} in \eqref{Grto0},
to obtain an identity without $\tilde g.$  Using the limit
\eqref{limghat} and that $\hat g$ and $g$ are regular, we obtain
\begin{eqnarray*}
 \int\limits_{\partial \Omega}\frac{\partial \tilde G }{\partial n}\frac{\partial G }{\partial
 n}n\,ds=\lim\limits_{r\to 0}\frac 1{r^{N-1}}\int\limits_{\partial B_r}\frac{1}{\sigma_N}\nabla
 \hat
g\,ds=\nabla \hat g(0,0)=\nabla \tilde \phi(0),
\end{eqnarray*}
where the last equality follows by observation after Theorem
\ref{tmain1}.
\end{proof}


\section{Proof of Lemma \ref{mainlemma}}

Let us recall the problem \eqref{emusystem1}--\eqref{emusystem2},
\begin{eqnarray}
-\Delta u_{\e,\mu}&=& v_{\e,\mu}^p\quad\mbox{in}\quad \Omega_\e \\
-\Delta v_{\e,\mu}&=& u_{\e,\mu}^{q_\e}\quad\mbox{in}\quad \Omega_\e \\
u_{\e,\mu}&=& v_{\e,\mu}=0\quad \mbox{on}\quad \partial \Omega_\e
\end{eqnarray}
where $\Omega_\e=(\Omega-x_\e)/\mu^{1-\e/2}.$
Let $\bar R>0.$ We define $\sigma(p):=2+N-p(N-2),$ and the scalar
function
$$
J(|y|):=\begin{cases} 1 \quad \mbox{if}\quad  \sigma(p)<2,\cr |\log(|y|/\bar R)|\quad \mbox{if}\quad  \sigma(p)=2, \cr
|y|^{2-\sigma(p)}\quad \mbox{if}\quad \sigma(p)>2.
\end{cases}
$$
Note that $\sigma(p)\in[0,N)$ for $p\in (2/(N-2),(N+2)/(N-2)]$ and
$\sigma(q)\leq 0.$ We consider the transformations
$$
z_\e(y)=|y|^{2-N}v_{\e,\mu}\left(\frac y{|y|^2}\right)\quad
\mbox{and}\quad
w_\e(y)=\frac{|y|^{2-N}}{J(|y|)}u_{\e,\mu}\left(\frac
y{|y|^2}\right)
$$
in $\Omega_\e^*$, the image of $\Omega_\e$ under $x\mapsto x/|x|^2.$

The next lemma is equivalent to Lemma \ref{mainlemma}, using the asymptotic behaviour  \eqref{hulshof2}.

\begin{lemma} Let $(w_\e,z_\e)$ solving
\begin{eqnarray}
-\Delta  J(|y|) w_\e &=&|y|^{-\sigma(p)} z_\e^p \quad \mbox{in}\quad\Omega^*_\e\\
-\Delta z_\e&=& |y|^{-\sigma(q)+(q_\e-q)(N-2)}[J(|y|) w_\e]^{q_\e} \mbox{in}\quad\Omega^*_\e\\
w_\e=z_\e&=&0\quad \mbox{on}\quad\partial\Omega^*_\e.
\end{eqnarray}
Then for any fixed $R\in(0,\bar R),$ we have
$$
\|w_\e\|_{L^{\infty}(\Omega_\e^R)}+\|z_\e\|_{L^{\infty}(\Omega_\e^R)}\leq C
$$
where $\Omega_\e^R=\Omega_\e^*\cap B_R$, and $C=C(R)$ independent of $\e>0$ provided $\e$ is sufficiently small.
\end{lemma}

\begin{proof}
Given $R>0,$ let $w_0$ and $z_0$ be solutions of
$$
\Delta J(|y|)w_0=0\quad \mbox{in}\quad \Omega^R_\e\quad \mbox{and}\quad w_0=0,\:\mbox{on}\quad \partial \Omega_\e^*\quad
w_0=w_\e\quad \mbox{on}\quad \partial B_R,
$$
and
$$
\Delta z_0=0\quad \mbox{in}\quad \Omega^R_\e\quad \mbox{and}\quad
z_0=0,\:\mbox{on}\quad\partial \Omega^*_\e\quad z_0=z_\e\quad
\mbox{on}\quad \partial B_R.
$$
By the convergence in compact sets of $w_\e$ and $z_\e$ , see
\eqref{convcontinuos}, we have $|z_\e|+|\nabla
z_\e|+|w_\e|+|\nabla w_\e|\leq C$ in $|y|=R$ for $C$ independent
of $\e$. Therefore by the maximum principle, we get
$$
 |Jw_0|+|\nabla (Jw_0)|+|z_0|+|\nabla z_0|\leq C \quad \mbox{in}\quad \Omega_\e^R.
$$
Define $\tilde w=w_\e-w_0$ and $\tilde z=z_\e-z_0.$ We now write
\begin{eqnarray}\label{tildea}
-\Delta  J(|y|) \tilde w  &=& a(y) z_\e \quad
\mbox{in}\quad\Omega^R_\e\\ \label{tildeb}
-\Delta \tilde z&=&b(y) J(|y|) w_\e\quad  \mbox{in}\quad\Omega^R_\e\\
\tilde w &=& \tilde z \:=\: 0\quad
\mbox{on}\quad\partial\Omega^R_\e \label{tildebc}
\end{eqnarray}
where $a(y)=|y|^{-\sigma(p)}z_\e^{p-1}$ and $b(y)=|y|^{-\sigma(q)+(q_\e-q)(N-2)}[J(|y|) w_\e]^{q_\e-1}.$
Clearly by the maximum
principle $\tilde w\geq 0$ and $\tilde z\geq 0.$

Let $P(y)=a(y)$ and
$$
Q(y)=
\begin{cases}
\frac 1{M}b(y) \quad\mbox{for}\quad y\in B_R\setminus \bar B_r  \cr
b(y) \quad\mbox{for}\quad B_r
\end{cases}
$$
where $r\in (0,R)$ and $M>1$ both independent of $\e$ and to be
determined later. Then
$$
b(y)J(|y|)w_\e=Q(y)J(|y|)w_\e+f(y)
$$
where
$$
f(y)=(b(y)-Q(y))J(|y|)w_\e=
\begin{cases}
0 \quad\mbox{for}\quad y\in \Omega_\e\cap B_r \cr
(1-\frac 1{M})b(y)J(|y|)w_\e\quad\mbox{for}\quad y\in B_R\setminus \bar B_r
\end{cases}
$$
It is clear that $f\in L^\infty(\Omega_\e^R)$, in fact
$\|f\|_{L^\infty(\Omega_\e^R)}\leq (1-1/M)r^{-(2+N)}$ by using
that $w_\e(y)\leq Cr^{\sigma(p)-N}$ for $|y|\geq r$, when
$p<N/(N-2),$ and $w_\e(y)\leq Cr^{2-N}$ for $|y|\geq r$ when
$p>N/(N-2).$ A similar bound is obtained for $p=N/(N-2).$
Then we transform \eqref{tildea}--\eqref{tildeb} in the system
\begin{eqnarray*}
-\Delta J\tilde w&= &P z_\e \quad \mbox{in}\quad\Omega^R_\e \\
-\Delta \tilde z&= &Q J w_\e+f \quad \mbox{in}\quad\Omega^R_\e
\end{eqnarray*}
We define $\eta_2(y)=\chi_{w_\e\leq 2\tilde w}(y)$ and
$\eta_1(y)=\chi_{z_\e\leq 2\tilde z}(y)$ for $y\in \Omega_\e^R,$
we find
\begin{eqnarray*}
-\Delta J \tilde w&\leq &2\eta_1 P\tilde z+ f_1 \quad \mbox{in}\quad\Omega^R_\e\\
-\Delta \tilde z&\leq &2\eta_2 Q J\tilde w+ f_2 \quad \mbox{in}\quad\Omega^R_\e
\end{eqnarray*}
Where $f_1=(1-\eta_1)Pz_\e=\chi_{z_\e\leq 2 z_0}Pz_\e\leq 2Pz_0$
and $f_2=f +(1-\eta_2)QJw_\e$ where $(1-\eta_2)QJw_\e\leq 2QJw_0.$
We write the system in the form
\begin{eqnarray}
-\Delta J\tilde w&\leq  & 2\eta_1 P |y|^{\gamma} |y|^{-\gamma}\tilde z+ f_1 \quad \mbox{in}\quad\Omega^R_\e,\\
-|y|^{-\gamma}\Delta \tilde z&\leq  & 2\eta_2 Q |y|^{-\gamma} J
\tilde w+f_2|y|^{-\gamma} \quad \mbox{in}\quad\Omega^R_\e, \\
\tilde w &=& \tilde z \:=\: 0\quad
\mbox{on}\quad\partial\Omega^R_\e.
\end{eqnarray}
Let $u(y)\mapsto  2\eta_2 Q |y|^{-\gamma}u(y)$ and $u(y)\mapsto
2\eta_1 P|y|^{\gamma}u(y)$  be the multiplication operators
$\mathcal{ P}$ and $\mathcal{Q}$ respectively. Note that a
multiplication operator $\mathcal C$ with corresponding function
$c(y)\in L^s(\Omega^R_\e)$ is bounded from $L^{s_1}(\Omega^R_\e)$
to $L^{s_2}(\Omega^R_\e)$ with $1/s_2=1/s_1+1/s.$

Formally we define $-L$ the operator as $u(y)\mapsto
-|y|^{-\gamma}\Delta(|y|^{\gamma}u(y)).$ More precisely, in the
appendix, we define $(-\Delta)^{-1}$ and $(-L)^{-1},$ which by the
Hardy-Littlewood-Sobolev inequality are bounded, independently of
$\e,$ from $L^{m_1}(\Omega^R_\e)$ to $L^{m_2}(\Omega^R_\e)$ with
$1/m_1=1/m_2+2/N.$ Note that the image of these operators is a
function with zero-Dirichlet boundary condition, so they are
positive. Then we can write
$$
J\tilde w\leq (-\Delta )^{-1}{\mathcal P}(-L)^{-1}({\mathcal Q} (J
\tilde w)+|y|^{-\gamma}f_2)+(-\Delta)^{-1} f_1.
$$

Denoting by $K=(-\Delta )^{-1}{\mathcal P}(-L)^{-1}{\mathcal Q}$
and $h=K|y|^{-\gamma} f_2+(-\Delta)^{-1} f_1$ we have
$$
(I-K)J\tilde w\leq h
$$
The proof is complete finding $m$ large enough such that $h\in L^m(\Omega^R_\e)$ and $(I-K)$ is invertible from
$L^m(\Omega^R_\e)$ to $L^m(\Omega^R_\e).$

We can estimate $Q(y)|y|^{-\gamma}$ in
$L^{\frac{q+1}{q-1}}(\Omega_\e^R),$ for
$\gamma=2\sigma(p)/(p+1)\geq 0,$ and note that
$\gamma=-\sigma(q)/(q+1)$ using the Sobolev Hyperbola. Since
$v_{\e,\mu}\to V$ in $L^{q+1}(\R^N),$ we have
$$
\int\limits_{\Omega_\e^*}
[J(|y|)w_\e(y)-V(y/|y|^2)|y|^{2-N}]^{q+1}|y|^{-\sigma(q)} \,d y
\to 0  \quad \mbox{as}\quad \e\to 0.
$$
Therefore for any $\lambda, $ we can take $r$ small such that
$$
\int\limits_{\Omega_\e^r} [J w_\e]^{(q+1)\frac{q_\e-1}{q-1}}(y)
|y|^{-\sigma(q)} \,dy \leq \int\limits_{\Omega_\e^r} [J
w_\e]^{(q+1)}(y) |y|^{-\sigma(q)} \,dy
\leq\frac{\lambda}{2C(\delta)}
$$
and $M$ large such that for all $\e\leq \e_0$ we have
\begin{eqnarray}\nonumber
\int\limits_{\Omega_\e^R} [Q(s)|y|^{-\gamma}]^{\frac{q+1}{q-1}}\,dy&\leq& C(\delta)\int\limits_{\Omega_\e^r} [J w_\e]^{(q+1)\frac{q_\e-1}{q-1}}|y|^{-\sigma(q)}\,dy\\
 &&+
 \frac{C(\delta)}{M^\frac{q+1}{q-1}}\int\limits_{B_R\setminus B_r} [J w_\e]^{(q+1)\frac{q_\e-1}{q-1}}|y|^{-\sigma(q)}\,dy \leq \lambda
\label{cotalambda}
\end{eqnarray}
where we have used $b(y)\leq C(\delta)[J w_\e]^{q_\e-1}$ with
$\delta$ given by  Lemma \ref{deltalemma}.

Now we show that $K$ is bounded from  $L^m(\Omega^R_\e)$ to $L^m(\Omega^R_\e).$
\begin{eqnarray*}
\|KJ\tilde w\|_{L^m(\Omega_\e^R)}&\leq& C_1\|\mathcal
P(-L)^{-1}\mathcal Q J \tilde
w\|_{L^r(\Omega^R_\e)}\qquad\qquad\qquad\qquad\qquad
\\
&\leq&  C_1\||y|^{\gamma}2\eta_1
P\|_{L^{\frac{p+1}{p-1}}(\Omega_\e^R)}\|(-L)^{-1}\mathcal QJ\tilde
w\|_{L^{r'}(\Omega^R_\e)}  \\
&\leq & C_1 \||y|^{\gamma} 2\eta_1
P\|_{L^{\frac{p+1}{p-1}}(\Omega_\e^R)}C_2\|\mathcal Q J\tilde
w\|_{L^{s'}(\Omega_\e^R)} \\ &\leq & C_1C_2\||y|^{\gamma} 2\eta_1
P\|_{L^{\frac{p+1}{p-1}}(\Omega_\e^R)}\||y|^{-\gamma}2\eta_2
Q\|_{L^{\frac{q+1}{q-1}}(\Omega_\e^R)}\|J\tilde
w\|_{L^{m'}(\Omega_\e^R)} \\
&\leq & \overline{C} \||y|^{\gamma}
P\|_{L^{\frac{p+1}{p-1}}(\Omega_\e^R)}\||y|^{-\gamma}
Q\|_{L^{\frac{q+1}{q-1}}(\Omega_\e^R)}\|J\tilde
w\|_{L^{m'}(\Omega_\e^R)}
\end{eqnarray*}
with $\frac 1{r}=\frac 1{m}+\frac 2{N},$ so $r'>1$ implies
$m>N/(N-2).$  $\frac 1{r}=\frac {p-1}{p+1}+\frac 1{r'}$ and $\frac
1{s'}=\frac {1}{r'}+\frac 2{N},$ so condition b) in
\eqref{HLSconst} implies $N-2+N/m>2N/(p+1)$ and $s'>1$ implies
$m>(q+1)/2$ so a) in \eqref{HLSconst} holds since $\gamma>0$ and
$\frac 1{s'}=\frac {q-1}{q+1}+\frac 1{m'}.$ Since
$$
\frac{q-1}{q+1}+\frac {p-1}{p+1}=\frac 4{N},\quad \mbox{we have}
\quad m'=m.
$$
By
$$
\int\limits_{\Omega_\e^*}
[z_\e(y)-U(y/|y|^2)|y|^{2-N}]^{p+1}|y|^{-\sigma(p)} \,d y
\to 0  \quad \mbox{as}\quad \e\to 0,
$$
we deduce that $ \||y|^{\gamma-\sigma(p)}z_\e^{p-1}
\|_{L^{\frac{p+1}{p-1}}(\Omega_\e^R)}=\||y|^\gamma
P\|_{L^{\frac{p+1}{p-1}}(\Omega_\e^R)}\leq C(\e_0) $ with
$C(\e_0)>0$ and for all $\e\in(0,\e_0).$ Since $\lambda$ in
\eqref{cotalambda} can be arbitrarely small then the norm of $K$
is small and so $I-K\colon L^m(\Omega^R_\e)\to L^m(\Omega^R_\e)$
invertible for $m$ large.  We have that
$$
\||y|^{-\gamma}f_2\|_{L^{m}(\Omega_\e^R)}\leq r^{-\gamma}\|f_2\|_{L^{\infty}(\Omega_r^R)}({\rm meas}(\Omega_r^R))^{1/m}
$$
is bounded, since $f_2$ is zero outside $\Omega_r^R$ and
$$
\|\Delta^{-1}f_1\|_{L^{m}(\Omega_\e^R)}\leq C_1\| f_2\|_{L^{r}(\Omega_\e^R)}\leq \| f_1\|_{L^{\infty}(\Omega_r^R)}({\rm
meas}(\Omega_\mu^R))^{1/r}\leq C(z_0)({\rm meas}(B_R))^{1/r}
$$

This implies $\|J w\|_{L^m(\Omega_\e^R)}\leq M $ for every $m$ large, and consequently for every $m\geq 1.$ (Use the $w_0$ to get
that $\|J w_\e\|_{L^m(\Omega_\e^R)}\leq M $). Now we have that
$$
-\Delta \tilde z = b(y)J w_\e= |y|^{-\sigma(q)+(q-q_\e)(N-2)}[J w_\e]^{q_\e}.
$$
Since $\sigma(q)\leq 0,$ if we take $m$ large such that $m q_\e>N/2$ then
\begin{equation}\label{boundzeta}
\|\tilde z\|_{L^\infty(\Omega^R_\e)}\leq \tilde M\quad\mbox{and
therefore}\quad \|z_\e\|_{L^\infty(\Omega^R_\e)}\leq M
\end{equation} for some $M$ independent of $\e.$ We study now each
case of $J$ separately. We have
\begin{equation}\label{boundzeta2}
-\Delta Jw_\e=|y|^{-\sigma(p)}z_\e^p\quad \mbox{in}\quad
\Omega_\e^*.
\end{equation}
a) In the case $J=1,$ since $\sigma(p)<2,$ using
\eqref{boundzeta}, we have $-\Delta\tilde w_\e\in L^{q}(\Omega)$
for any $q\in (N/2, N/\sigma(p)).$ By regularity, we get
$$
\|w_\e\|_{L^\infty(\Omega^R_\e)}\leq M.
$$
b) For $J(|y|)=-\log(|y|/\bar R)>\log(\bar R/R),$ we have
$$
-\Delta \tilde w -\frac{\nabla J}{J}\nabla \tilde w -\frac{\Delta J}{J}\tilde w=\frac 1{J|y|^{2}}z_\e^p\quad \mbox{in}\quad
\Omega^R_\e
$$
or equivalently
$$
-\Delta \tilde w +\frac{1}{J|y|^2}(y,\nabla \tilde w) +\frac{1}{J|y|^2}(N-2)\tilde w=\frac 1{J|y|^{2}}z_\e^p\quad \mbox{in}\quad
\Omega^R_\e.
$$
Using \eqref{boundzeta}, we can take $u=\tilde w-M$ with
$M=\sup\limits_{\e>0}\sup\limits_{y\in \Omega^R_\e}
z_\e^p(y)/(N-2),$ and we get
$$
-J|y|^2\Delta u +(y,\nabla u) +(N-2)u \leq 0\quad \mbox{in}\quad \Omega^R_\e.
$$
Since $u=-M< 0$ on the boundary,  $ u\leq 0$ in $\Omega_\e^R.$ This gives $w_\e\leq M$  in $\Omega_\e^R.$

For the remaining case, $p<N/(N-2)$ we have
$$
-\Delta \tilde w -\frac{\nabla J}{J}\nabla \tilde w -\frac{\Delta J}{J}\tilde w=\frac 1{|y|^{2}}z_\e^p\quad \mbox{in}\quad
\Omega^R_\e.
$$
As before, defining $u=\tilde w-M$ with
$M=\sup\limits_{\e>0}\sup\limits_{y\in \Omega^R_\e}
z_\e^p/[(\sigma(p)-2)(N-\sigma(p))]$ then
$$
-|y|^2\Delta u - (2-\sigma(p))(y,\nabla u) - (2-\sigma(p))(N-\sigma(p)) u \leq 0\quad \mbox{in}\quad \Omega^R_\e
$$
Since $u=-M< 0$ on the boundary, $ u\leq 0$ in $\Omega_\e^R.$ This implies $w_\e\leq M$  in $\Omega_\e^R.$


\end{proof}

\section{Appendix}

Let $N>2.$ Let $h$ and $v$  be a function in
$L^{s'}(\Omega_\e^R).$ Given the Green's function $G$ solution of
$-\Delta G(x,\cdot) = \delta_{x}$ in $\Omega_\e^R,$ $G(x,\cdot)=0$
on $\partial \Omega_\e^R,$ we define
$$
(-\Delta)^{-1} h(\xi)=\int\limits_{\Omega_\e^R} G(x,\xi) h(x)\,dx
\quad \xi\in \Omega_\e^R.
$$
and
$$
(-L)^{-1} v(\xi)=|\xi|^{-\gamma}\int\limits_{\Omega_\e^R} G(x,\xi)
|x|^\gamma v(x)\,dx \quad \xi\in \Omega_\e^R.
$$
Note that $G$ is positive, so both operators are positive. We know
that  $(-\Delta)^{-1}$ is bounded, independently of $\e,$ from
$L^{s'}(\Omega^R_\e)$ to $L^{r'}(\Omega^R_\e)$ with
$1/r'=1/s'-2/N.$ Next we prove the same result for $(-L)^{-1}.$ By
the weighted Hardy-Littlewood-Sobolev inequality \cite{CLO,L}, for
$|\xi|^{-\gamma} f\in L^{s'}(\Omega_\e^R),$ we have that
$$
\|\xi^{-\gamma}(-\Delta)^{-1} f\|_{L^{r'}(\Omega_\e^R)}\leq
2\||\xi|^{-\gamma}\int\limits_{\Omega_\e^R} \frac
C{|x-\xi|^{N-2}}f(x)\,dx\|_{L^{r'}(\Omega_\e^R)}\leq C\|
|\xi|^{-\gamma} f\|_{L^{s'}(\Omega_\e^R)}
$$
for $1<s'<r'<\infty,$ with $1/r'=1/s'-2/N$ and
\begin{equation}\label{HLSconst}
a)\quad -\gamma<N(1-1/s')=N-2-N/r'\quad\mbox{and}\quad b)\quad \gamma<N/r'.
\end{equation}
In other words, for any $v\in L^{s'}(\Omega_\e^R),$ we have
\begin{eqnarray}\nonumber
\|(-L)^{-1}
v\|_{L^{r'}(\Omega_\e^R)}&=&\||\xi|^{-\gamma}(-\Delta)^{-1}
|x|^{\gamma} v\|_{L^{r'}(\Omega_\e^R)} \\
&\leq& \nonumber 2 \||\xi|^{-\gamma}\int\limits_{\Omega_\e^R}
\frac
C{|x-\xi|^{N-2}}|x|^{\gamma} v(x)\,dx\|_{L^{r'}(\Omega_\e^R)}\\
&\leq & C\| v \|_{L^{s'}(\Omega_\e^R)}.\label{L-1}
\end{eqnarray}

\begin{lemma}\label{Lpreg} Let $u$ solve
$$
\begin{cases} -\Delta u= f \quad \mbox{in}\quad \Omega \subset \R^N\cr
u=0\quad \mbox{on}\quad \partial \Omega
\end{cases}
$$
Let $\omega$ be a neighborhood of $\partial \Omega.$ Then
$$
\|u\|_{W^{1,q}(\Omega)}+\|\nabla u\|_{C^{0,\alpha}(\omega')}\leq
C(\|f\|_{L^1(\Omega)}+\|f\|_{L^\infty(\omega)})
$$
for $q<N/(N-1),$ $\alpha\in (0,1)$ and $\omega'\subset \omega$ is a strict subdomain of $\omega.$
\end{lemma}
\section*{Acknowledgements}

The author thank Prof. J. D\'avila for many useful discussions.

\end{document}